\documentclass[final]{dmtcs-episciences}

\usepackage[T1]{fontenc}
\usepackage[utf8]{inputenc}
\usepackage{microtype}
\usepackage{lmodern}
\usepackage{amsmath,amssymb,amsfonts,amsthm,amscd,enumerate,epsfig}
\usepackage{setspace}
\usepackage{mathrsfs}
\usepackage{fancyhdr} 
\usepackage{mathtools} 
\usepackage{listings}  
\usepackage{float} 
\usepackage[utf8]{inputenc}
\usepackage{array}
\usepackage[caption=false,font=normalsize,labelfont=sf,textfont=sf]{subfig}
\usepackage{textcomp}
\usepackage{stfloats}
\usepackage{url}
\usepackage{verbatim}
\usepackage{epsfig,fancybox}
\usepackage{algorithmicx}
\usepackage{algorithm}
\usepackage[noend]{algpseudocode}
\usepackage{hyperref}
\usepackage[round]{natbib}
\usepackage{tikz}
\usetikzlibrary{calc}
\usepackage[utf8]{inputenc}
\usepackage{graphicx}
\usepackage{textcomp}
\usepackage{multicol} 

\newcommand\DoubleLine[7][4pt]{
\path(#2)--(#3)coordinate[at start](h1)coordinate[at end](h2);
\draw[#4]($(h1)!#1!90:(h2)$)-- node [auto=left] {#5} ($(h2)!#1!-90:(h1)$);
\draw[#6]($(h1)!#1!-90:(h2)$)-- node [auto=right] {#7} ($(h2)!#1!90:(h1)$);}

\usepackage[round]{natbib}

\author{
    Juan Carlos Garc\'ia-Altamirano\thanks{Supported by Conacyt Grant 963921.} 
    \and Mika Olsen\thanks{Supported by Conacyt Project 47510664} 
    \and Jorge Cervantes-Ojeda}
\affiliation{
  UAM CUAJIMALPA, M\'exico}
\publicationdata{vol. 25:1}{2023}{3}{10.46298/dmtcs.10189}{2022-10-21; 2022-10-21; 2023-01-31}{2023-02-03}

\title[Symmetric cycle of length 5 using Haj\'os construction and adapted Rank Genetic Algorithm]{How to construct the symmetric cycle of length 5 using Haj\'os construction with an adapted Rank Genetic Algorithm}

\begin{document}
\maketitle

\begin{abstract}
In 2020 Bang-Jensen et al. generalized the Haj\'os join of two graphs to the class of digraphs and generalized several results for vertex colorings in digraphs. 
Although, as a consequence of these results, a digraph can be obtained by Haj\'os constructions (directed Haj\'os join and identifying non-adjacent vertices), determining the Haj\'os constructions to obtain the digraph is a complex problem. In particular, Bang-Jensen et al. posed the problem of determining the Haj\'os operations to construct the symmetric 5-cycle from the complete symmetric digraph of order 3 using only Haj\'os constructions.
We successfully adapted a rank-based genetic algorithm to solve this problem by the introduction of innovative recombination and mutation operators from graph theory. The Haj\'os Join became the recombination operator and the identification of independent vertices became the mutation operator. In this way, we were able to obtain a sequence of only 16 Haj\'os operations to construct the symmetric cycle of order 5.   
\end{abstract}

\section{Introduction}
A genetic algorithm can solve mathematical problems that are intractable by exhaustive search, as is the case of the problem described here.
Evolutionary algorithms have been applied to a wide variety of engineering problems, for instance: \cite{Dasgupta2013EvolutionaryAI}, and they have also been applied to mathematics problems: \cite{Jong1989UsingGA} showed that Genetic Algorithms (GA) can be used to solve NP-complete problems, \cite{JAKOBS1996165} applied a GA to a geometry problem, \cite{POURRAJABIAN201311483} solved nonlinear algebraic equations by using a GA, and \cite{Cervantes2019} applied a Rank Genetic Algorithm (Rank GA)  to the graph theory problem of finding the rainbow connection number of a graph.   \cite{https://doi.org/10.1002/jgt.20003} proved that for a given digraph $D$, deciding whether the dichromatic number of a digraph is at most 2 is NP-Complete. Therefore, our motivation is to try using the Rank GA heuristic on this problem. To our best knowledge, this problem has not been approached in this way nor using any other heuristic algorithms.

In a genetic algorithm, we have an initial population of individuals. In our case, each individual represents a digraph obtained by  Haj\'os constructions of the complete symmetric digraph on 3 vertices $D(K_3)$. Each individual in the population is evaluated according to a fitness function. This function measures the ability of the individual for reaching a predetermined objective. Once individuals are evaluated, the next generation is obtained by applying genetic operators to the population inspired by the evolution in nature, such as cross-over between individuals, selection of the fittest individuals, and mutations. This procedure is repeated until the genetic algorithm finds an individual that achieves the required objective or until a maximum number of generations is reached.

A vertex coloring of a digraph is acyclic if there are no monochromatic directed cycles. The dichromatic number of a digraph $D$ was introduced in \cite{NEUMANNLARA1982265} as the minimum number of colors of an acyclic coloring of $D$, denoted by $dc(D)$. The dichromatic number of a digraph is an extension of the chromatic number of a graph and several concepts and results for the chromatic number of a graph have been extended to digraphs using the dichromatic number. For instance, perfect digraphs by \cite{Andres2015PerfectD},  dichromatic polynomial by \cite{10.1007/s00373-022-02484-0}, Brooks theorem for digraphs by \cite{Harutyunyan2011StrengthenedBT}, bounds in terms of the girth of a graph by \cite{CorderoMichel2018NewBF}, flow theory  by \cite{HOCHSTATTLER2017160}, and diachromatic number by \cite{AraujoPardo2017TheDN}.
A digraph $D$ is {\bf r-critical} if $dc(D)=r$ and $dc(H)<k$ for every proper subdigraph $H$ of $D$. 
 \cite{BangJensen2019HajsAO} extended the well-known Haj\'os construction for graphs to digraphs: any $r$-critical digraph can be obtained by Haj\'os constructions using complete symmetric digraphs on $r$ vertices, $D(K_r)$. Although the result was proved, it is not a trivial task to obtain even simple digraphs such as symmetric cycles of odd length. In particular, the authors left as an open problem how to construct the symmetric cycle $D(C_5)$ using directed Haj\'os joins and identifying non-adjacent vertices.
Another interesting problem is constructing digraphs of minimum order for a given dichromatic number $r$. The 3-dichromatic tournaments of minimum order were characterized in \cite{NEUMANNLARA1994233}.

\subsection{Definitions}
We consider finite digraphs without loops and without multiple arcs. For all definitions not given here we refer the reader to the book of Chartrand, Lesniak and Zhang \cite{Chartrand2010GraphsD}.  
Let $D$ be a digraph with vertex set $V (D)$ and arc set $A(D)$. The {\bf in-neighborhood} of a vertex $u$ is $N^-(u) = \{v\in V(D)\mid vu\in A(D)\}$  and the {\bf out-neighborhood} of a vertex $u$ is $N^+(u) = \{v\in V(D)\mid uv\in A(D)\}$. 
Two vertices in a digraph $D$ are {\bf independent} if there are no arcs between them in $D$, a set of vertices $X$ is independent in a digraph $D$ if any pair of vertices of $X$ is independent in $D$.
An arc $uv\in A(D)$ is {\bf symmetric} ({\bf asymmetric}) if $vu\in A(D)$ ($vu\notin A(D)$), and a digraph is symmetric (bidirected graph) if every arc of $D$ is a symmetric arc.
The {\bf symmetric digraph} $D(G)$, of the graph $G$, is the digraph obtained by replacing each edge with a symmetric arc.     

The Haj\'os construction was defined for digraphs in 2020 by Bang-Jensen et al. \cite{BangJensen2019HajsAO} as an extension of the well-known Haj\'os construction \cite{Hajos,https://doi.org/10.1002/(SICI)1097-0118(199901)30:1<37::AID-JGT5>3.0.CO;2-V,https://doi.org/10.1002/(SICI)1097-0118(199712)26:4<211::AID-JGT5>3.0.CO;2-T} for graphs. The class of {\bf Haj\'os-k-constructible digraphs} is defined as the smallest family of digraphs that contains all complete digraphs of order $k$ and is closed under directed Haj\'os join and identifying independent vertices.    

The digraph $H$ obtained by {\bf identifying} a non-empty set $I$ of independent vertices is defined as the digraph $H=D-I$ adding a new vertex $v$  and adding all arcs from $v$ to $N^+_D(I)=\bigcup\limits_{u\in I} N^+_D(u)$ and all arcs from $N^-_D(I)=\bigcup\limits_{u\in I} N^-_D(u)$ to $v$. The new vertex $v$ may preserve the label of one of the vertices of the independent set $I$.

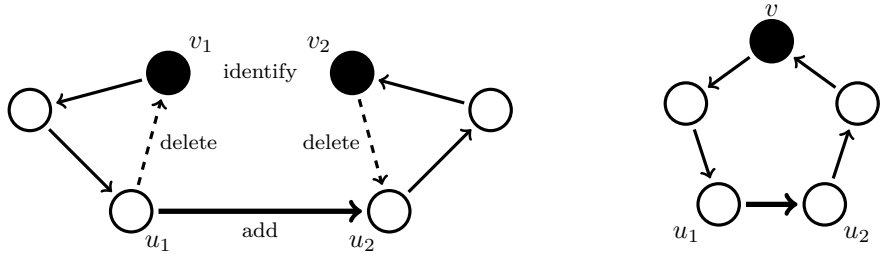
\begin{figure*}[ht]
\centering
\beginpgfgraphicnamed{example}
\begin{center}
\begin{tikzpicture}[myn/.style={circle,very thick,draw,inner sep=0.11cm,outer sep=2.5pt},label distance=-1.5mm]]

\node [myn,circle,label=45:$v_1$,fill=black,minimum size=0.55cm] 
(0) at (45.000000:1.1000) {};

\node [myn,circle,minimum size=0.55cm] 
(1) at (165.000000:1.1000) {};

\node [myn,circle,label=285:$u_1$,minimum size=0.55cm] 
(2) at (285.000000:1.1000) {};

\draw[very thick,black][->] (0) -- (1);
\draw[very thick,black][->] (1) -- (2);
\draw[very thick,black][dashed,->] (2) -- (0) node [pos=0.5,right,font=\footnotesize] {delete};

\pgftransformshift{\pgfpoint{4.00000cm}0cm}

\node [myn,circle,label=135:$v_2$,fill=black,minimum size=0.55cm] 
(3) at (135.000000:1.1000) {};

\node [myn,circle,label=255:$u_2$,minimum size=0.55cm] 
(4) at (255.000000:1.1000) {};

\node [myn,circle,minimum size=0.55cm] 
(5) at (15.000000:1.1000) {};

\draw[very thick,black][dashed,->] (3) -- (4) node [pos=0.5,left,font=\footnotesize] {delete};
\draw[very thick,black][->] (4) -- (5);
\draw[very thick,black][->] (5) -- (3);

\draw[line width = 2.0pt,black][->] (2) -- (4) node [pos=0.5,below,font=\footnotesize] {add};
\draw[line width = 0.0pt,white][-] (0) -- (3) node [pos=0.5,font=\footnotesize,black] {identify};

\pgftransformshift{\pgfpoint{4.80000cm}0cm}

\node [myn,circle,label=90:$v$,fill=black,minimum size=0.55cm] 
(0) at (90.000000:1.2000) {};

\node [myn,circle,minimum size=0.55cm] 
(1) at (162.000000:1.2000) {};

\node [myn,circle,label=234:$u_1$,minimum size=0.55cm] 
(2) at (234.000000:1.2000) {};

\node [myn,circle,label=306:$u_2$,minimum size=0.55cm] 
(3) at (306.000000:1.2000) {};

\node [myn,circle,minimum size=0.55cm] 
(4) at (378.000000:1.2000) {};

\draw[very thick,black][->] (0) -- (1);
\draw[very thick,black][->] (1) -- (2);
\draw[very thick,black][->, line width = 2.0pt] (2) -- (3);
\draw[very thick,black][->] (3) -- (4);
\draw[very thick,black][->] (4) -- (0);
\end{tikzpicture}
\end{center}
\endpgfgraphicnamed

\caption{The first digraph indicates the directed Haj\'os join of two directed triangles and the second is the resulting digraph.}
\label{Fig_Hajos}
\end{figure*}

We use Figure \ref{Fig_Hajos} to illustrate the definition of directed Haj\'os Join. Let $D_1$ and $D_2$ be two disjoint digraphs. Let $u_1v_1\in A(D_1)$ and $v_2u_2\in A(D_2)$. The {\bf directed Haj\'os join} $D=(D_1,u_1,v_1)\triangledown(D_2,v_2,u_2)$ or, briefly $D=D_1\triangledown D_2$ of $D_1$ and $D_2$ is defined as the disjoint union of $D_1$ and $D_2$ and deleting both arcs $u_1v_1$ and $v_2u_2$, identifying the vertices $v_1$ and $v_2$ to a new vertex $v$ and adding the arc $u_1u_2$. The vertex $v$ may be denoted by $v_1$, $v_2$ or $v$. 

\section{The problem}
\cite{BangJensen2019HajsAO} proved that any $k$-critical digraph is Haj\'os-$k$-constructible. In the same paper they posed the following question (Question 19): 

{How can a bidirected (symmetric) $C_5$ ($D(C_5)$) be constructed from copies of $D(K_3)$ by only using directed Haj\'os joins and identifying non-adjacent vertices?}

To answer this question, we use an evolutionary algorithm that exclusively applies directed Haj\'os joins and identifications of non-adjacent vertices. The algorithm starts with a population that contains only copies of $D(K_3)$ and is guided using a fitness function whose optimum is in the digraph $D(C_5)$, which is the symmetric cycle of length 5. It is important to stress that as a consequence of Theorem 4 in \cite{BangJensen2019HajsAO}, we already know that the symmetric cycle of order 5 can be obtained from $D(K_3)$ using only Haj\'os operations, but we do not know how. From the algorithm, we save the operations that are being carried out by the evolutionary algorithm to obtain each individual, in order to reconstruct the operations needed to obtain the optimum, which in this case is the digraph $D(C_5)$. 

Observe that the digraphs in the populations obtained after applying Haj\'os constructions, may have different order and size. 

\section{Using a Genetic Algorithm}
In this section, we describe how we use a genetic algorithm to solve our problem.
\subsection{The Rank GA}
The Rank GA is an algorithm that was first presented in  \cite{Cervantes2009} as a solution to the limitations of existing mutation rate heuristics. It has been successfully applied to solve several problems, for instance  \cite{Cervantes2019, jaimes2018simultaneous, 5585961, 6119011}. The claim is that this algorithm has a very good balance between exploration and exploitation by applying all genetic operations to the population based (or depending) on the fitness rank of each individual.  

In the Rank GA, the individuals of the population are ranked from best to worst in terms of their fitness before the genetic operators (selection, recombination, and mutation) are applied. The application of these genetic operators depends on the rank of each individual in the population. The top ranked individuals tend to vary less than the bottom ranked ones. This is to make the latter try to escape from the local optima of the fitness function. Also, top ranked individuals tend to be cloned more than others who tend to disappear.

\subsection{Adaptation to our problem}
We use an adapted version of the Rank GA by \cite{Cervantes2009} to solve the problem. The adapted algorithm (Haj\'os RankGA) pseudocode is given in Algorithm \ref{alg:rankGA} which we describe below. In that algorithm, the population size is equal to  $50$ and the number of generations has no limit. Here we give an overview of the features taken from \cite{Cervantes2009} and the modifications and adaptions made to that. For details, please consult that original paper.

\begin{algorithm}[H]
\caption{Haj\'os Rank GA}\label{alg:rankGA}
\begin{algorithmic}
\Procedure{rankGA()}{}
\For{each individual $i$}
   \State   individual$_i \leftarrow D(K_3)$
\EndFor
\State   Evaluation and Sort
\While{not end Criteria met}
   \State   rankRecombination()
   \State   Evaluation and Sort
   \State   rankMutation()
   \State   Evaluation and Sort
   \State   rankSelection() and Sort
\EndWhile
\EndProcedure
\end{algorithmic}
\end{algorithm}

\begin{algorithm}[H]
\caption{rank Recombination}\label{alg:rankRecombination}
\begin{algorithmic} 
\Procedure{rankRecombination()}{}\
\For{$i$ in [0..individuals.size()-1] step 2} 
   \State   arc1 = random(arcs of individual$_i$)
   \State   arc2 = random(arcs of individual$_{i+1}$)
   \State   individuals.add(
    Haj\'{o}s($i, (i+1)$  mod individuals.size(), arc1, arc2))
\EndFor
\EndProcedure
\end{algorithmic}
\end{algorithm}

\begin{algorithm}[H]
\caption{Haj\'{o}s operation}
\label{alg:hajosOperation}
\begin{algorithmic}
\State   
\Procedure{Haj\'{o}s(int $i$, int $j$, arc1, arc2)}{}\
\State   result $\leftarrow$ union(individual$_i$, individual$_j$) \Comment{disjoint union} 
\State   result.deleteArcs(arc1, arc2)
\State   identify(result, arc2.v2, arc1.v1)
\State   result.addArc(arc1.u1, arc2.u2)
\EndProcedure
\end{algorithmic}
\end{algorithm}

\begin{algorithm}[H]
\caption{identification operation}
\label{alg:identificationOperation}
\begin{algorithmic}
\State
\Procedure{identify(individual, $u$, $v$)}{}\
\State   individual.vertex($u$).addArcs(arcs of $v$)
\State   individual.delete($v$)
\EndProcedure
\end{algorithmic}
\end{algorithm}

\begin{algorithm}[H]
\caption{rank Mutation}
\label{alg:rankMutation}
\begin{algorithmic}
\State
\Procedure{rankMutation()}{}\
\For{$i$ in [0..2individuals.size()-1] } 
    \State   $r \leftarrow i$/(2individuals.size-1)
    \State   nonAdjacentPairs $\leftarrow$ binomial numVertices$_i$,2)-(asim$_i$+digon$_i$)
    \State   numAttempts $\leftarrow$ (r)nonAdjacentPairs 
    \For{$j$ in [1..numAttempts]}
          \State    select 1 random vertex $u$ from individual$_i$
          \If {$u$ has non adjacent vertices} 
             \State   $v$ = random(non adjacent vertices of $u$)
             \State   identify(individual$_i$, $v$, $u$)
          \EndIf
    \EndFor
\EndFor
\EndProcedure
\end{algorithmic}
\end{algorithm}

\begin{algorithm}[ht]
\caption{rank Selection}
\label{alg:rankSelection}
\begin{algorithmic}
\State   
\Procedure{rankSelection()}{}\
\State   clones $\leftarrow$ null
\For{$i$ in [0..2individuals.size()-1]} 
    \State   $r \leftarrow i$/(2individuals.size()-1)
    \State   $n \leftarrow \lfloor P(1-r)^{(2P-1)} \rfloor$
    \For{$j$ in [0..$n-1$] }
       \State   clones.add(individual$_i$)  \Comment {individual $i$ is cloned $n$ times}
    \EndFor
\EndFor
\State   $i \leftarrow 0$
\While {clones.size() < individuals.size()} 
    \State   $r \leftarrow i$/(2individuals.size()-1)
    \State   $n \leftarrow P(1-r)^{(2P-1)}$ 
    \State   $f \leftarrow n - \lfloor n \rfloor$  \Comment {$f$ is the fractional part of $n$}
    \If {random(0,1) $< f$}
       \State   clones.add(individual$_i$)  \Comment {one extra clone of individual $i$}
    \EndIf
    \State   $i \leftarrow (i+1)$ mod 2individuals.size()
\EndWhile
\State   individuals $\leftarrow$ clones  \Comment {replace population}
\EndProcedure
\end{algorithmic}
\end{algorithm}

\subsubsection{Representation}
Each individual in the population represents a digraph $D$ by its adjacency matrix $A$, defined as follows

Let $D$ be a digraph of order $n$ with $V(D) = \{v_0,v_1,\dots,v_{n-1}\}$ being its set of vertices.  The {\bf adjacency matrix} of $D$ is the matrix $A_{n\times n}$, where $A[i,j] = 1$ if $v_iv_j$ is an  arc of the digraph $D$ and $A[i,j] = 0$ otherwise.

\subsubsection{Fitness Function}
The fitness function to be minimized that we use is:
\begin{equation}\label{fitness}
   ft = |n-5|+2\frac{a}{n}+|n-s|+15T_S+5T,
\end{equation}
where $n$ is the order of the individual, $a$ is the number of its asymmetric arcs, $s$ is the number of its symmetric arcs, $T_S$ is the number of its $D(K_3)$ subdigraphs and $T$ is the number of its triangles with at least one non-symmetric arc. 

The term $|n-5|$ is introduced to favor those digraphs of order 5. 
The term $2\frac{a}{n}$ is introduced to favor digraphs with a low density of asymmetric arcs relative to their order. The term $|n-s|$ favors digraphs with the same number of vertices and symmetric arcs. 
The terms $15T_S$ and $5T$ are introduced to favor those digraphs with no triangles (symmetric and non-symmetric). 
Observe that if the fitness function is equal to 0, the digraph has 5 vertices, no asymmetric arcs, 5 symmetric arcs, and no triangles. The unique 3-dichromatic digraph that satisfies these conditions is the symmetric cycle of length 5, $D(C_5)$. 

\subsubsection{Recombination}

Recombination is done based on the rank of the individuals (see Algorithm  \ref{alg:rankRecombination}), so the population needs to be sorted by fitness value before any recombination.
Mating is done in such a way that the best individual is recombined only with the second best. The third best is recombined only with the fourth best, etc.

We use the directed Haj\'os join between two digraphs as our  {Recombination} operator.

Procedure $Haj\acute{o}s$ in Algorithm \ref{alg:hajosOperation} takes two individuals and one arc from each individual and returns the Haj\'os join between them. Procedure \emph{identify} in Algorithm \ref{alg:identificationOperation} takes an individual and a pair of independent vertices and produces another individual that is the input individual with the identification operation done between the given vertices. 

\subsubsection{Mutation}
The amount of mutation individuals get is a function of their rank (see Algorithm  \ref{alg:rankMutation}). The particular function to be used takes into account that the best individual should have 0 mutations and the worst individual should have the maximum number of mutations. First, we calculate the number of non-adjacent pairs of vertices and then we take a fraction of this to determine the number of mutation attempts to be done on the individual. 
This function is:

\begin{equation}
\label{eq:numAttempts}
numAttempts = (r)numNonAdjacentPairs
\end{equation}where 
\begin{equation}
\label{eq:rank}
r=\frac{i}{2popSize-1}
\end{equation}where $i$ is the index of the individual in the sorted population (starting at 0). Here we use $2popSize$ because after recombination we have doubled the population.
A single mutation attempt selects one vertex randomly and, if it has non-adjacent vertices, one of those is selected randomly, or else the attempt is aborted. Once a pair of non-adjacent vertices was selected, the mutation is performed as the identification of these two (see Algorithm \ref{alg:identificationOperation}). 

\subsubsection{Selection}
Rank selection defines that the number of clones of an individual is a function of its rank  (see Algorithm \ref{alg:rankSelection}). This function is monotonously decreasing and it is expected to produce a population of half the size of the previous population.
The number of clones for $individual_i$ is:
\begin{equation}
clones=P(1-r)^{2P-1}
\end{equation}
where $r$ is given by equation \ref{eq:rank} and $P$ is the selective pressure parameter that in our case is $P=3$.

The integer part of $numClones$ determines a minimum number of clones for that individual, whereas the fractional part of it determines the probability that this individual gets one extra clone. The probabilistic part is performed after the integer part and only until the required population size is reached.

\subsubsection{Stored information}
Each population is stored and each individual is assigned its {\bf origin} that describes how the individual is obtained in terms of the previous population. In case the individual was not modified, the origin indicates from which individual of the previous population it comes from, and in case the individual was modified, the origin indicates how it was obtained from one or two individuals from the previous population. Due to this information, we can reconstruct the operations that lead to the solution of the problem.

\section{Results}
The Haj\'os Rank GA used almost 5000 generations to obtain the symmetric 5-cycle, but we do not need 5000 Haj\'os operations. 
In each generation, the RankGA \ref{alg:rankGA} stores the population and the operations from which each individual was obtained, with this information we can recursively reconstruct the steps used by the algorithm to generate a particular individual.
That is, we can reconstruct recursively each individual from $D(K_3)$ no matter the number of generations.  
In what follows, we describe the 16 Haj\'os operations we obtained from the algorithm to construct the symmetric 5-cycle from $D(K_3)$ using Haj\'os operations. 
Thus, we answer the question posed in \cite{BangJensen2019HajsAO}. 

\begin{figure*}[ht]
\centering
\beginpgfgraphicnamed{example}
\begin{center}
\begin{tikzpicture}[myn/.style={circle,very thick,draw,inner sep=0.11cm,outer sep=2.5pt},label distance=-1.5mm]]

\node [myn,circle,label=45:$v_0$,fill=black,minimum size=0.55cm] 
(0) at (45.000000:1.1000) {};

\node [myn,circle,label=165:$v_1$,minimum size=0.55cm] 
(1) at (165.000000:1.1000) {};

\node [myn,circle,label=285:$v_2$,minimum size=0.55cm] 
(2) at (285.000000:1.1000) {};

\DoubleLine{0}{1}{<-,very thick}{}{->,very thick}{}
\DoubleLine{1}{2}{<-,very thick}{}{->,very thick}{}
\DoubleLine{2}{0}{<-,very thick}{}{dashed,->,very thick}{};

\pgftransformshift{\pgfpoint{3.50000cm}0cm}

\node [myn,circle,label=135:$v'_0$,fill=black,minimum size=0.55cm] 
(3) at (135.000000:1.1000) {};

\node [myn,circle,label=255:$v'_1$,minimum size=0.55cm] 
(4) at (255.000000:1.1000) {};

\node [myn,circle,label=15:$v'_2$,minimum size=0.55cm] 
(5) at (15.000000:1.1000) {};

\DoubleLine{3}{5}{<-,very thick}{}{->,very thick}{}
\DoubleLine{4}{5}{<-,very thick}{}{->,very thick}{}
\DoubleLine{3}{4}{<-,very thick}{}{dashed,->,very thick}{}
\draw[line width = 2.0pt,black][->] (2) -- (4) node [pos=0.5,below,font=\footnotesize] {};

\pgftransformshift{\pgfpoint{4.800000cm}0cm}

\node [myn,circle,label=90:$v_0$,fill=black,minimum size=0.55cm] 
(0) at (90.000000:1.2000) {};

\node [myn,circle,label=162:$v_1$,minimum size=0.55cm] 
(1) at (162.000000:1.2000) {};

\node [myn,circle,label=234:$v_2$,minimum size=0.55cm] 
(2) at (234.000000:1.2000) {};

\node [myn,circle,label=306:$v_3$,minimum size=0.55cm] 
(3) at (306.000000:1.2000) {};

\node [myn,circle,label=378:$v_4$,minimum size=0.55cm] 
(4) at (378.000000:1.2000) {};

\DoubleLine{0}{1}{<-,very thick}{}{->,very thick}{}
\draw[very thick,black][->] (3) -- (0);
\DoubleLine{0}{4}{<-,very thick}{}{->,very thick}{}
\DoubleLine{1}{2}{<-,very thick}{}{->,very thick}{}
\draw[very thick,black][->, line width = 2.0pt] (2) -- (3);
\DoubleLine{3}{4}{<-,very thick}{}{->,very thick}{}
\draw[very thick][->] (0) -- (2);

\end{tikzpicture}
\end{center}
\endpgfgraphicnamed

\caption{The first digraph indicates the construction of the digraph $D_0=(H,v_2,v_0)\triangledown(H',v'_0,v'_1)$ 
using two copies of $D(K_3)$, and the second is the resulting digraph.}
\label{Fig_step_1}
\end{figure*}
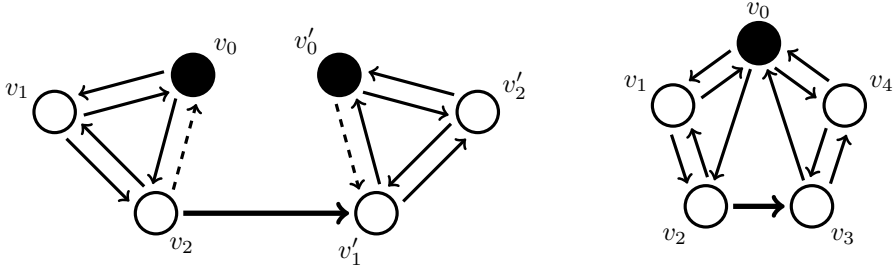

The following sequence of directed Haj\'os joins and vertex identifications is the one obtained by the adapted version of the Haj\'os Rank GA: \bigskip

The initial step is the directed Haj\'os join between two disjoint copies $H$ and $H'$ of $D(K_3)$, where $V(H)=\{v_0,v_1,v_2\}$ and $V(H')=\{v'_0,v'_1,v'_2\}$. In the left digraph in Figure \ref{Fig_step_1}, we indicate the arcs that must be deleted with a dotted line, and we indicate the vertices that we identify with a black background, and the right digraph in Figure \ref{Fig_step_1} is the resulting digraph where $v_0$ is the vertex obtained by identifying the vertices $v_0$ and $v_0'$, and with the vertices that correspond to the vertices of $H'$ relabeled in cyclic order as $v_3$ and $v_4$.  

\begin{figure*}[ht] 
\centering
{
\begin{tikzpicture}[myn/.style={circle,very thick,draw,inner sep=0.11cm,outer sep=2.5pt},label distance=-1.5mm]]

\node [myn,circle,label=90:$v_0$,fill=black,minimum size=0.55cm] 
(0) at (90.000000:1.2000) {};

\node [myn,circle,label=162:$v_1$,fill=black!60,minimum size=0.55cm] (1) at (162.000000:1.2000) {};

\node [myn,circle,label=234:$v_2$,fill=black!35,minimum size=0.55cm] (2) at (234.000000:1.2000) {};

\node [myn,circle,label=306:$v_3$,fill=black!15,minimum size=0.55cm] (3) at (306.000000:1.2000) {};

\node [myn,circle,label=378:$v_4$,minimum size=0.55cm] 
(4) at (378.000000:1.2000) {};

\DoubleLine{0}{1}{<-,very thick}{}{->,very thick}{}
\DoubleLine{0}{4}{<-,very thick}{}{->,very thick}{}
\DoubleLine{1}{2}{<-,very thick}{}{->,very thick}{}
\draw[very thick,black][->] (4) -- (1);
\DoubleLine{2}{3}{<-,very thick}{}{->,very thick}{}
\DoubleLine{3}{4}{<-,very thick}{}{->,very thick}{}
\draw[very thick][dashed,->] (3) -- (0);

\pgftransformshift{\pgfpoint{3.900000cm}0cm}

\node [myn,circle,label=90:$v'_0$,fill=black!60,minimum size=0.55cm] (5) at (90.000000:1.2000) {};

\node [myn,circle,label=162:$v'_1$,fill=black!35,minimum size=0.55cm] (6) at (162.000000:1.2000) {};

\node [myn,circle,label=234:$v'_2$,fill=black!15,minimum size=0.55cm] (7) at (234.000000:1.2000) {};

\node [myn,circle,label=306:$v'_3$,minimum size=0.55cm] 
(8) at (306.000000:1.2000) {};

\node [myn,circle,label=378:$v'_4$,fill=black,minimum size=0.55cm] (9) at (378.000000:1.2000) {};

\DoubleLine{5}{6}{<-,very thick}{}{->,very thick}{}
\draw[very thick,black][->] (8) -- (5);
\DoubleLine{5}{9}{<-,very thick}{}{->,very thick}{}
\DoubleLine{6}{7}{<-,very thick}{}{->,very thick}{}
\DoubleLine{7}{8}{<-,very thick}{}{->,very thick}{}
\DoubleLine{8}{9}{<-,very thick}{}{->,very thick}{}
\draw[very thick][dashed,->] (9) -- (6);
\draw[line width = 2.0pt,black][->] (3) -- (6);

\pgftransformshift{\pgfpoint{4.50000cm}0cm}

\node [myn,circle,label=90:$v_0$,fill=black,minimum size=0.55cm] 
(0) at (90.000000:1.2000) {};

\node [myn,circle,label=162:$v_1$,fill=black!60,minimum size=0.55cm] (1) at (162.000000:1.2000) {};

\node [myn,circle,label=234:$v_2$,fill=black!35,minimum size=0.55cm] (2) at (234.000000:1.2000) {};

\node [myn,circle,label=306:$v_3$,fill=black!15,minimum size=0.55cm] (3) at (306.000000:1.2000) {};

\node [myn,circle,label=378:$v_4$,minimum size=0.55cm] 
(4) at (378.000000:1.2000) {};

\DoubleLine{0}{1}{<-,very thick}{}{->,very thick}{}
\DoubleLine{0}{4}{<-,very thick}{}{->,very thick}{}
\DoubleLine{1}{2}{<-,very thick}{}{->,very thick}{}
\draw[very thick,black][->] (4) -- (1);
\DoubleLine{3}{4}{<-,very thick}{}{->,very thick}{}
\pgftransformshift{\pgfpoint{0.000000cm}0cm}
\DoubleLine{3}{2}{<-,very thick}{}{->,line width = 2.0pt}{}

\end{tikzpicture}
}

\caption{The first digraph indicates the construction of the digraph $D_2=(D_1,v_3,v_0)\triangledown(D'_1,v'_4,v'_1)$ using two copies of  $D_1$, and the second is the resulting digraph.}
\label{Fig_step_3}
\end{figure*}
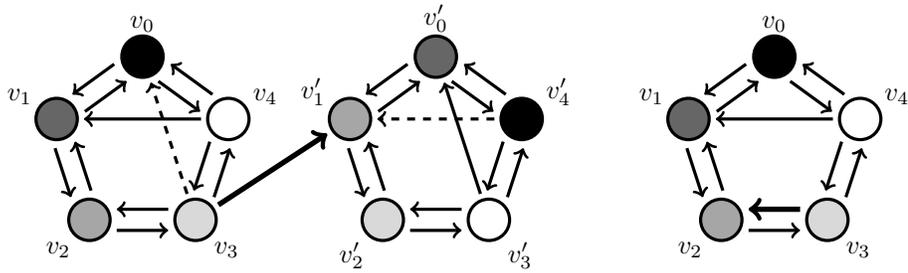

\begin{figure*}[h!] 
\centering
{
\begin{tikzpicture}[myn/.style={circle,very thick,draw,inner sep=0.11cm,outer sep=2.5pt},label distance=-1.5mm]]

\node [myn,circle,label=90:$v_0$,minimum size=0.55cm] 
(0) at (90.000000:1.2000) {};

\node [myn,circle,label=162:$v_1$,fill=black,minimum size=0.55cm] (1) at (162.000000:1.2000) {};

\node [myn,circle,label=234:$v_2$,fill=black!60,minimum size=0.55cm] (2) at (234.000000:1.2000) {};

\node [myn,circle,label=306:$v_3$,fill=black!35,minimum size=0.55cm] (3) at (306.000000:1.2000) {};

\node [myn,circle,label=378:$v_4$,fill=black!15,minimum size=0.55cm] 
(4) at (378.000000:1.2000) {};

\DoubleLine{0}{1}{<-,very thick}{}{->,very thick}{}
\DoubleLine{0}{4}{<-,very thick}{}{->,very thick}{}
\DoubleLine{1}{2}{<-,very thick}{}{->,very thick}{}
\DoubleLine{2}{3}{<-,very thick}{}{->,very thick}{}
\DoubleLine{3}{4}{<-,very thick}{}{->,very thick}{}
\draw[very thick][dashed,->] (4) -- (1);

\pgftransformshift{\pgfpoint{3.900000cm}0cm}

\node [myn,circle,label=90:$v'_0$,fill=black!60,minimum size=0.55cm] (5) at (90.000000:1.2000) {};

\node [myn,circle,label=162:$v'_1$,fill=black!35,minimum size=0.55cm] (6) at (162.000000:1.2000) {};

\node [myn,circle,label=234:$v'_2$,fill=black!15,minimum size=0.55cm] (7) at (234.000000:1.2000) {};

\node [myn,circle,label=306:$v'_3$,minimum size=0.55cm] 
(8) at (306.000000:1.2000) {};

\node [myn,circle,label=378:$v'_4$,fill=black,minimum size=0.55cm] (9) at (378.000000:1.2000) {};

\DoubleLine{5}{6}{<-,very thick}{}{->,very thick}{}
\DoubleLine{5}{9}{<-,very thick}{}{->,very thick}{}
\DoubleLine{6}{7}{<-,very thick}{}{->,very thick}{}
\DoubleLine{7}{8}{<-,very thick}{}{->,very thick}{}
\DoubleLine{8}{9}{<-,very thick}{}{->,very thick}{}
\draw[very thick][dashed,->] (9) -- (6);
\draw[line width = 2.0pt,black][->] (4) -- (6);

\pgftransformshift{\pgfpoint{4.50000cm}0cm}

\node [myn,circle,label=90:$v_0$,minimum size=0.55cm] 
(0) at (90.000000:1.2000) {};

\node [myn,circle,label=162:$v_1$,fill=black,minimum size=0.55cm] (1) at (162.000000:1.2000) {};

\node [myn,circle,label=234:$v_2$,fill=black!60,minimum size=0.55cm] (2) at (234.000000:1.2000) {};

\node [myn,circle,label=306:$v_3$,fill=black!35,minimum size=0.55cm] (3) at (306.000000:1.2000) {};

\node [myn,circle,label=378:$v_4$,fill=black!15,minimum size=0.55cm] 
(4) at (378.000000:1.2000) {};

\DoubleLine{0}{1}{<-,very thick}{}{->,very thick}{}
\DoubleLine{0}{4}{<-,very thick}{}{->,very thick}{}
\DoubleLine{1}{2}{<-,very thick}{}{->,very thick}{}
\DoubleLine{2}{3}{<-,very thick}{}{->,very thick}{}
\pgftransformshift{\pgfpoint{0.000000cm}0cm}
\DoubleLine{4}{3}{<-,very thick}{}{->,line width = 2.0pt}{}

\end{tikzpicture}
}
\caption{The first digraph indicates the construction of the digraph $D_3=(D_2,v_4,v_1)\triangledown(D'_2,v'_4,v'_1)$ using two copies of  $D_2$, and the second is the resulting digraph.}
\label{Fig_step_4}
\end{figure*}
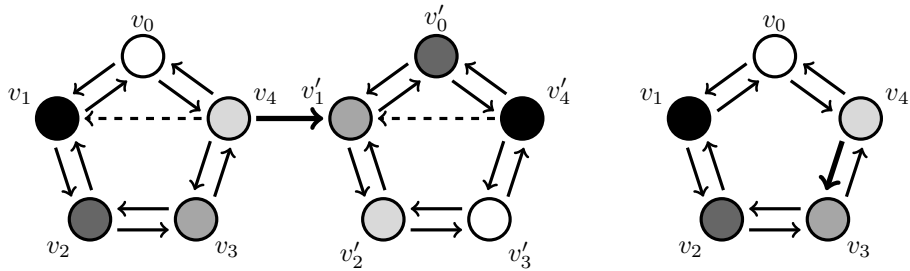

In each of the following steps, we consider 5 Haj\'os operations: a directed Haj\'os join and four vertex identifications. In Figures \ref{Fig_step_1} to  \ref{Fig_step_4}, the left digraph indicates the directed Haj\'os join, where the dotted arcs must be deleted, the thick arrow added and the black vertices identified, and the shades of gray indicate the four pairwise vertex identifications. The digraph on the right indicates the result of these operations.

The operation in Figure \ref{Fig_step_2} uses as input the result in Figure \ref{Fig_step_1}, the operation in Figure \ref{Fig_step_3} does so with the result in Figure \ref{Fig_step_2} and the one in Figure \ref{Fig_step_4} does so with the result in Figure \ref{Fig_step_3}.

\begin{figure*}[ht] 
\centering
\begin{tikzpicture}[myn/.style={circle,very thick,draw,inner sep=0.11cm,outer sep=2.5pt},label distance=-1.5mm]]

\node [myn,circle,label=90:$v_0$,fill=black!15,minimum size=0.55cm] 
(0) at (90.000000:1.2000) {};

\node [myn,circle,label=162:$v_1$,minimum size=0.55cm] 
(1) at (162.000000:1.2000) {};

\node [myn,circle,label=234:$v_2$,fill=black,minimum size=0.55cm] 
(2) at (234.000000:1.2000) {};

\node [myn,circle,label=306:$v_3$,fill=black!60,minimum size=0.55cm] (3) at (306.000000:1.2000) {};

\node [myn,circle,label=378:$v_4$,fill=black!35,minimum size=0.55cm] (4) at (378.000000:1.2000) {};

\DoubleLine{0}{1}{<-,very thick}{}{->,very thick}{}
\draw[very thick,black][->] (3) -- (0);
\DoubleLine{0}{4}{<-,very thick}{}{->,very thick}{}
\DoubleLine{1}{2}{<-,very thick}{}{->,very thick}{}
\draw[very thick,black][->] (2) -- (3);
\DoubleLine{3}{4}{<-,very thick}{}{->,very thick}{}
\draw[very thick][dashed,->] (0) -- (2);

\pgftransformshift{\pgfpoint{3.900000cm}0cm}

\node [myn,circle,label=90:$v'_0$,fill=black!35,minimum size=0.55cm] (5) at (90.000000:1.2000) {};

\node [myn,circle,label=162:$v'_1$,fill=black!15,minimum size=0.55cm] (6) at (162.000000:1.2000) {};

\node [myn,circle,label=234:$v'_2$,minimum size=0.55cm] 
(7) at (234.000000:1.2000) {};

\node [myn,circle,label=306:$v'_3$,fill=black,minimum size=0.55cm] 
(8) at (306.000000:1.2000) {};

\node [myn,circle,label=378:$v'_4$,fill=black!60,minimum size=0.55cm] (9) at (378.000000:1.2000) {};

\DoubleLine{5}{6}{<-,very thick}{}{->,very thick}{}
\draw[very thick,black][->] (5) -- (7);
\DoubleLine{5}{9}{<-,very thick}{}{->,very thick}{}
\DoubleLine{6}{7}{<-,very thick}{}{->,very thick}{}
\draw[very thick,black][->] (7) -- (8);
\DoubleLine{8}{9}{<-,very thick}{}{->,very thick}{}
\draw[very thick][dashed,->] (8) -- (5);
\draw[line width = 2.0pt,black][->] (0) -- (5);

\pgftransformshift{\pgfpoint{4.500000cm}0cm}

\node [myn,circle,label=90:$v_0$,fill=black!15,minimum size=0.55cm] 
(0) at (90.000000:1.2000) {};

\node [myn,circle,label=162:$v_1$,minimum size=0.55cm] 
(1) at (162.000000:1.2000) {};

\node [myn,circle,label=234:$v_2$,fill=black,minimum size=0.55cm] 
(2) at (234.000000:1.2000) {};

\node [myn,circle,label=306:$v_3$,fill=black!60,minimum size=0.55cm] (3) at (306.000000:1.2000) {};

\node [myn,circle,label=378:$v_4$,fill=black!35,minimum size=0.55cm] (4) at (378.000000:1.2000) {};

\DoubleLine{0}{1}{<-,very thick}{}{->,very thick}{}
\draw[very thick,black][->] (3) -- (0);
\DoubleLine{1}{2}{<-,very thick}{}{->,very thick}{}
\draw[very thick][->] (4) -- (1);
\DoubleLine{2}{3}{<-,very thick}{}{->,very thick}{}
\DoubleLine{3}{4}{<-,very thick}{}{->,very thick}{}
\DoubleLine{0}{4}{<-,very thick}{}{->,line width = 2.0pt}{}

\end{tikzpicture}

\caption{The first digraph indicates the construction of the digraph $D_1=(D_0,v_0,v_2)\triangledown(D'_0,v'_3,v'_0)$ using two copies of  $D_0$, and the second is the resulting digraph.}
\label{Fig_step_2}
\end{figure*}
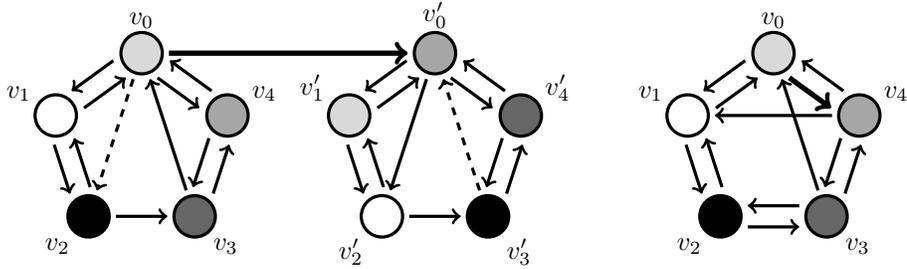

Observe that the resulting digraph in Figures \ref{Fig_step_1} and \ref{Fig_step_2} is a symmetric 5-cycle with two diagonals, the resulting digraph in Figure \ref{Fig_step_3} is a symmetric 5-cycle with only one diagonal and the resulting digraph in Figure \ref{Fig_step_4} is a symmetric 5-cycle without diagonals.

\section{Conclusions}
Graph and digraph theory is a field of application of GAs that has been little explored, although there are examples where GAs have been applied successfully, for instance \cite{BOUAZZI20214337,Cervantes2019}. 
The sequence of operations obtained by the GA inspired the authors to develop a method to obtain any symmetric odd cycle. The description of the operations of that method is out of the scope of this article, but it will be described in a forthcoming article by the authors. 
Thus, this article reiterates the viability of using genetic algorithms, with innovative ways for their adaptation, as an important heuristic in the development of graph theory.

\section*{Acknowledgment.} 
 Research supported by CONACYT Project 47510664 and the Department of Applied Mathematics and Systems (DMAS) at Metropolitan Autonomous University (UAM) Cuajimalpa, and 
Juan Carlos Garc\'ia-Altamirano was supported by Conacyt Grant 963921.

\nocite{*}
\bibliographystyle{abbrvnat}
\bibliography{Hajos}

\begin{thebibliography}{24}
\providecommand{\natexlab}[1]{#1}
\providecommand{\url}[1]{\texttt{#1}}
\expandafter\ifx\csname urlstyle\endcsname\relax
  \providecommand{\doi}[1]{doi: #1}\else
  \providecommand{\doi}{doi: \begingroup \urlstyle{rm}\Url}\fi

\bibitem[Andres and Hochst{\"a}ttler(2015)]{Andres2015PerfectD}
S.~D. Andres and W.~Hochst{\"a}ttler.
\newblock Perfect digraphs.
\newblock \emph{Journal of Graph Theory}, 79:\penalty0 21--29, 2015.

\bibitem[Araujo-Pardo et~al.(2017)Araujo-Pardo, Montellano-Ballesteros, Olsen,
  and Rubio-Montiel]{AraujoPardo2017TheDN}
G.~Araujo-Pardo, J.~J. Montellano-Ballesteros, M.~Olsen, and C.~Rubio-Montiel.
\newblock The diachromatic number of digraphs.
\newblock \emph{Electron. J. Comb.}, 25:\penalty0 \#P3.51, 2017.

\bibitem[Bang-Jensen et~al.(2019)Bang-Jensen, Bellitto, Schweser, and
  Stiebitz]{BangJensen2019HajsAO}
J.~Bang-Jensen, T.~Bellitto, T.~Schweser, and M.~Stiebitz.
\newblock Haj{\'o}s and ore constructions for digraphs.
\newblock \emph{Electron. J. Comb.}, 27:\penalty0 \#P1.63, 2019.

\bibitem[Bokal et~al.(2004)Bokal, Fijavz, Juvan, Kayll, and
  Mohar]{https://doi.org/10.1002/jgt.20003}
D.~Bokal, G.~Fijavz, M.~Juvan, P.~M. Kayll, and B.~Mohar.
\newblock The circular chromatic number of a digraph.
\newblock \emph{Journal of Graph Theory}, 46\penalty0 (3):\penalty0 227--240,
  2004.

\bibitem[Bouazzi et~al.(2021)Bouazzi, Hammami, and Bouamama]{BOUAZZI20214337}
K.~Bouazzi, M.~Hammami, and S.~Bouamama.
\newblock Application of an improved genetic algorithm to hamiltonian circuit
  problem.
\newblock \emph{Procedia Computer Science}, 192:\penalty0 4337--4347, 2021.

\bibitem[Cervantes and Stephens(2009)]{Cervantes2009}
J.~Cervantes and C.~R. Stephens.
\newblock Limitations of existing mutation rate heuristics and how a rank ga
  overcomes them.
\newblock \emph{IEEE Transactions on Evolutionary Computation}, 13:\penalty0
  369--397, 2009.

\bibitem[Cervantes et~al.(2010)Cervantes, Sánchez, and González]{5585961}
J.~Cervantes, M.~Sánchez, and P.~P. González.
\newblock Emerging traits in the application of an evolutionary algorithm to a
  scalable bioinformatics problem.
\newblock \emph{IEEE Congress on Evolutionary Computation}, pages 1--8, 2010.

\bibitem[Cervantes-Ojeda et~al.(2019)Cervantes-Ojeda, del Carmen
  G{\'o}mez~Fuentes, Gonz{\'a}lez-Moreno, and Olsen]{Cervantes2019}
J.~Cervantes-Ojeda, M.~del Carmen G{\'o}mez~Fuentes, D.~Gonz{\'a}lez-Moreno,
  and M.~Olsen.
\newblock Rainbow connectivity using a rank genetic algorithm: Moore cages with
  girth six.
\newblock \emph{J. Appl. Math.}, 2019:\penalty0 4073905:1--4073905:7, 2019.

\bibitem[Chartrand et~al.(2015)Chartrand, Lesniak, and
  Zhang]{Chartrand2010GraphsD}
G.~Chartrand, L.~M. Lesniak, and P.~Zhang.
\newblock \emph{Graphs \& Digraphs, Sixth Edition}.
\newblock Chapman and Hall/CRC., 2015.

\bibitem[Cordero-Michel and Galeana-S{\'a}nchez(2018)]{CorderoMichel2018NewBF}
N.~Cordero-Michel and H.~Galeana-S{\'a}nchez.
\newblock New bounds for the dichromatic number of a digraph.
\newblock \emph{Discret. Math. Theor. Comput. Sci.}, 21, 2018.

\bibitem[Dasgupta and Michalewicz(2013)]{Dasgupta2013EvolutionaryAI}
D.~Dasgupta and Z.~Michalewicz.
\newblock \emph{Evolutionary Algorithms in Engineering Applications}.
\newblock Springer, 2013.

\bibitem[Flores and Cervantes(2011)]{6119011}
D.~Flores and J.~Cervantes.
\newblock Rank based evolution of real parameters on noisy fitness functions:
  Evolving a robot neurocontroller.
\newblock \emph{2011 10th Mexican International Conference on Artificial
  Intelligence}, pages 72--76, 2011.

\bibitem[Gonz\'{a}lez-Moreno et~al.(2022)Gonz\'{a}lez-Moreno,
  Hern\'{a}ndez-Ortiz, Llano, and Olsen]{10.1007/s00373-022-02484-0}
D.~Gonz\'{a}lez-Moreno, R.~Hern\'{a}ndez-Ortiz, B.~Llano, and M.~Olsen.
\newblock The dichromatic polynomial of a digraph.
\newblock \emph{Graphs and Combinatorics}, 38\penalty0 (3):\penalty0 \#85,
  2022.

\bibitem[Haj\'os(1961)]{Hajos}
G.~Haj\'os.
\newblock \"uber eine konstruktion nicht n-f\"arbbarer graphen.
\newblock \emph{Wiss. Z. Martin Luther Univ. Halle-Wittenberg, Math. Natur.
  Reihe}, 10:\penalty0 116--117, 1961.

\bibitem[Harutyunyan and Mohar(2011)]{Harutyunyan2011StrengthenedBT}
A.~Harutyunyan and B.~Mohar.
\newblock Strengthened brooks' theorem for digraphs of girth at least three.
\newblock \emph{Electron. J. Comb.}, 18, 2011.

\bibitem[Hochstättler(2017)]{HOCHSTATTLER2017160}
W.~Hochstättler.
\newblock A flow theory for the dichromatic number.
\newblock \emph{European Journal of Combinatorics}, 66:\penalty0 160--167,
  2017.
\newblock Selected papers of EuroComb15.

\bibitem[Jakobs(1996)]{JAKOBS1996165}
S.~Jakobs.
\newblock On genetic algorithms for the packing of polygons.
\newblock \emph{European Journal of Operational Research}, 88\penalty0
  (1):\penalty0 165--181, 1996.

\bibitem[Jensen and
  Royle(1999)]{https://doi.org/10.1002/(SICI)1097-0118(199901)30:1<37::AID-JGT5>3.0.CO;2-V}
T.~R. Jensen and G.~F. Royle.
\newblock Hajós constructions of critical graphs.
\newblock \emph{Journal of Graph Theory}, 30\penalty0 (1):\penalty0 37--50,
  1999.

\bibitem[Jong and Spears(1989)]{Jong1989UsingGA}
K.~A.~D. Jong and W.~M. Spears.
\newblock \emph{Using Genetic Algorithms to Solve NP-Complete Problems}.
\newblock IEEE, 1989.

\bibitem[L\'opez-Jaimes et~al.(2018)L\'opez-Jaimes, Cervantes-Ojeda,
  G{\'o}mez-Fuentes, and Alvarado-Gonz{\'a}lez]{jaimes2018simultaneous}
A.~L\'opez-Jaimes, J.~Cervantes-Ojeda, M.~C. G{\'o}mez-Fuentes, and A.~M.
  Alvarado-Gonz{\'a}lez.
\newblock Simultaneous evolution of neuro-controller for multiple car-like
  robots.
\newblock \emph{Research in Computer Science}, 147\penalty0 (10):\penalty0
  29--44, 2018.

\bibitem[Neumann-Lara(1982)]{NEUMANNLARA1982265}
V.~Neumann-Lara.
\newblock The dichromatic number of a digraph.
\newblock \emph{Journal of Combinatorial Theory, Series B}, 33\penalty0
  (3):\penalty0 265--270, 1982.

\bibitem[Neumann-Lara(1994)]{NEUMANNLARA1994233}
V.~Neumann-Lara.
\newblock The 3 and 4-dichromatic tournaments of minimum order.
\newblock \emph{Discrete Mathematics}, 135\penalty0 (1):\penalty0 233--243,
  1994.

\bibitem[Pourrajabian et~al.(2013)Pourrajabian, Ebrahimi, Mirzaei, and
  Shams]{POURRAJABIAN201311483}
A.~Pourrajabian, R.~Ebrahimi, M.~Mirzaei, and M.~Shams.
\newblock Applying genetic algorithms for solving nonlinear algebraic
  equations.
\newblock \emph{Applied Mathematics and Computation}, 219\penalty0
  (24):\penalty0 11483--11494, 2013.

\bibitem[Urquhart(1997)]{https://doi.org/10.1002/(SICI)1097-0118(199712)26:4<211::AID-JGT5>3.0.CO;2-T}
A.~Urquhart.
\newblock The graph constructions of hajós and ore.
\newblock \emph{Journal of Graph Theory}, 26\penalty0 (4):\penalty0 211--215,
  1997.

\end{thebibliography}
\label{sec:biblio}

\end{document}